\documentclass[11pt,a4paper]{article}
\usepackage{amsmath,amsfonts,amssymb,color}
\usepackage{amsthm}
\usepackage[hidelinks]{hyperref}
\usepackage{tikz}
\usepackage{pgfplots}
\usepackage{enumerate}
\usepackage[english]{babel}

\usepackage{makeidx}
\makeindex

\date{August 29, 2023}

\usepackage[utf8]{inputenc}
\usepackage[T1]{fontenc}

\author{Bruno Luiz Santos Correia and Marc Troyanov\footnote{Institut  de Math\'{e}matiques EPFL, 1015 Lausanne, Switzerland, 
bruno.santoscorreia@epfl.ch,  \  marc.troyanov@epfl.ch }}

\newtheorem{theorem}{\textrm{Theorem}}[section]
\newtheorem{lemma}[theorem]{\rm\bf Lemma}
\newtheorem{proposition}[theorem]{\rm\bf Proposition}
\newtheorem{corollary}[theorem]{\rm\bf Corollary}
\newtheorem{definition}[theorem]{\rm\bf Definition}
\newtheorem{example}[theorem]{\rm\bf Example}
\newtheorem{remark}[theorem]{\rm\bf Remark}
\newtheorem{question}[theorem]{\rm\bf Question}

\numberwithin{equation}{section}

\newcommand{\LL}{\mathcal{L}}
\newcommand{\U}{\mathcal{U}}

\newcommand{\Z}{\mathbb{Z}}
\newcommand{\N}{\mathbb{N}}
\newcommand{\R}{\mathbb{R}}

\DeclareMathOperator{\id}{id}
\DeclareMathOperator{\ee}{{e}}
 
\title{Isoperimetry in  Finitely Generated Groups}

\begin{document}
\maketitle

\begin{abstract}\parindent=0pt
We revisit the  isoperimetric inequalities for finitely generated groups  introduced and studied by N. Varopoulos, T. Coulhon and L. Saloff-Coste. Namely we show that  a lower bound on  the isoperimetric quotient of finite subsets  in a finitely generated group is given by the $\U-$transform of its growth function, which is a variant of the Legendre transform. From this lower bound, we obtain some asymptotic estimates for the F{\o}lner function of the group.  The paper also includes a discussion of some basic definitions from Geometric Group Theory and some basic properties of the $\U$-transform, including some computational techniques and its relation with the Legendre transform.
 
\medskip 

{\small  {\textit{Keywords: Finitely generated groups, growth, isoperimetric inequality, F{\o}lner function,  Lambert $W$-function.}}}

\smallskip 

{\small  {AMS Subject Classification: 20F65;  20F69, 05C63.}}
\end{abstract}

\tableofcontents

\newpage

\section{Introduction}

In the early 1990s, T. Coulhon  and  {L. Saloff-Coste} proved the following general isoperimetric  inequality relating the size of an arbitrary finite set $D$ in an infinite,  finitely generated group $\Gamma$ to the size of its boundary:
\begin{equation} \label{ineqCSC}
    \frac{|\partial_SD|}{|D|} \geq \frac{1}{4 |S|}\cdot \frac{1}{\phi_S (2 |D|)}.
\end{equation}
\index{Isoperimetric  inequality}
Here $S$ is a finite symmetric  set of generators for  $\Gamma$ and $\phi_S$ is the \textit{inverse growth function} associated with   the word metric in $\Gamma$ defined by $S$. Precise definitions are given in the next section. The remarkable feature of this inequality is that,  in some sense, the size of metric balls controls the  boundary size of arbitrary finite subsets of the group $\Gamma$. The proof of this inequality appeared in  \cite[Theorem 1]{CSC}, see also  \cite[Theorem 3.2]{PSC}; however, a version of this result for groups of polynomial growth has been proved earlier by N. Varopoulos in \cite{Varopoulos}.
%
M. Gromov  later improved the constant in the above inequality; he proved: 
\begin{equation} \label{ineqG}
  \frac{|\partial_SD|}{|D|} \geq \frac{1}{2}\cdot \frac{1}{\phi_S (2 |D|)},
\end{equation}
see Chapter 6, p. 346, in \cite{Gromov}. The proof also applies to finite groups, provided one assumes $|D| < \frac{1}{2} |\Gamma|$. 
 This version of the  isoperimetric inequality is now standard and appears in textbooks such as \cite[Theorem 14.95]{Ceccherini} and \cite[Theorem 6.29]{LyonsPeres}. 

\medskip

The constant $2$ appearing in the previous inequality is suboptimal; optimizing the argument in Gromov's proof, leads us to a  new formulation of the isoperimetric inequality in which   the inverse growth $\phi_S$ is replaced  by another transformation of  the growth function $\gamma$,  that we call the $\U$-\textit{transform} of that function.  In Theorem \ref{maintheoremA} we obtain the following inequality
$$
   \frac{|\partial_SD|}{|D|} \geq \U_{\gamma}(|D|),
$$
where  $\U_{\gamma}(t)$ is defined   in \eqref{defUT}. As explained in \S \ref{sec.legendre}, the   $\U$-transform  $\gamma \mapsto \U_{\gamma}$ is equivalent to the Legendre transform after some change of variables.  The main results of this paper  are Theorems   \ref{maintheoremA} and   \ref{maintheorem}, and various corollaries discussed in Section 3.

\medskip 

The rest of the paper is organized as follows: In Section 2 we recall some basic definitions from Geometric Group Theory. In Section 3  we formulate and prove the main result of the paper and discuss some of  its consequences. 
The basic properties of the $\U-$transform, including some computational techniques and its relation with the  Legendre transform are postponed till the last section of the paper.

\section{Preliminaries from Geometric Group Theory}

We recall here some   basic notions of geometric group theory and refer to the books \cite{Ceccherini}  and \cite{Harpe}  for a thorough  introduction to the subject. Starting with a group  $\Gamma$  generated by the finite symmetric set $S = S^{-1} \subset \Gamma$, we define: 

\bigskip

\textbf{(i) The word metric and the Cayley graph.}  
 For any $r\in \N$, we denote by $B_S(r) \subset \Gamma$ the subset of those elements in $\Gamma$ that can be written as a product of at most $r$ generators in $S$. The  right invariant \textit{word metric}  \index{Word metric}  on  $\Gamma$ with respect to the generating set $S$ is  the function $d_S : \Gamma \times \Gamma \to \mathbb{N}\cup \{0\}$ defined as
$$
  d_S(x,y) = \min \{r \in \mathbb{N}\cup \{0\}  \mid xy^{-1} \in B_S(r)\}.
$$
With this definition, $B_S(r)$ is the ball in the metric space $(\Gamma,d_S)$ centered at the identity element $e\in \Gamma$. 
We shall also denote by $\|x\|_S = d_S(e,x)$ the distance from $e$ to $x\in \Gamma$.

\smallskip

The \textit{Cayley graph} \index{Cayley graph}  of $\Gamma$ with respect to $S$ is the undirected graph $X(\Gamma,S)$ 
whose vertex set is the group $\Gamma$ itself and the edges are the pairs $\{x,y\} \subset \Gamma$ such that $d_S(x,y) = 1$ (i.e. such that $y=sx$ for some $s\in S$). Note that this is a regular graph of degree $d = |S|$. By construction, the distance $d_S$ between two  elements in the group $\Gamma$ is  the length of the shortest edge path joining them. 

\medskip 

As a simple example, consider  the Dihedral group  $D_4 = \langle r, s \mid r^4, s^2, (rs)^2  \rangle$. Note that $r^k s= sr^{4-k}$ for $0 \leq k \leq 4$,  and the corresponding Cayley graph is 
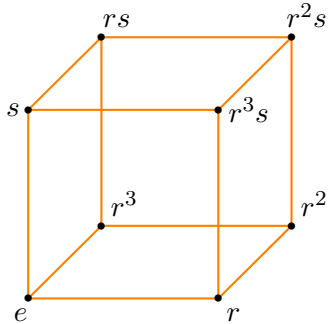
\begin{center}
\begin{tikzpicture}[scale=2.5]
  \coordinate (A) at (0,0,0);
  \coordinate (B) at (1,0,0);
  \coordinate (C) at (1,1,0);
  \coordinate (D) at (0,1,0);
  \coordinate (E) at (0,0,1);
  \coordinate (F) at (1,0,1);
  \coordinate (G) at (1,1,1);
  \coordinate (H) at (0,1,1);
  \foreach \x/\y in {A/B,B/C,C/D,D/A,E/F,F/G,G/H,H/E,A/E,B/F,C/G,D/H} {\draw[orange,thick]  (\x) -- (\y);
  }
\foreach \x in  {(0,0,0),(0,0,1),(1,0,0),(1,1,1),(0,1,0),(0,1,1),(1,1,0),(1,0,1)}{
    \node at \x {\tiny$\bullet$}; 
}
  \node[below,xshift=-0.1cm]  at (0,0,1) {$e$}; 
  \node[above,xshift=0.3cm]  at (0,0,0) {$r^3$}; 
  \node[above,xshift=0.3cm]  at (1,0,0) {$r^2$}; 
  \node[,below,xshift=0.2cm]  at (1,0,1) {$r$};
  \node[above,xshift=0.2cm]  at (0,1,0) {$rs$}; 
  \node[xshift=0.4cm]  at (1,1,1) {$r^3s$}; 
  \node[above,xshift=0.2cm]  at (1,1,0) {$r^2s$}; 
  \node[xshift=-0.2cm]  at (0,1,1) {$s$}; 
\end{tikzpicture}
\end{center}
Simple examples of infinite Cayley graphs are the infinite regular  tree of degree $2q$ for the free group $F_q$ on $q$ generators and the standard integer lattice in $\R^d$ for the free abelian group on $d$ generators.

\medskip

\textbf{(ii) The growth function.} The  \textit{growth function}   \index{Growth of a group}  of $\Gamma$   with respect to the generating set $S$  is the function   $\gamma_{_S}  : \N \to \N$ counting the number of elements in the ball $B(n)$. That is,  
$$
  \gamma_{_S} (n) =  |B_S(n)|,
$$
where we denote by $|D|$ the cardinality of a finite set $D$.  
Some examples of finitely generated groups and their growth are:
\begin{enumerate}[$\circ$]
\item A group has bounded growth  if and only if it is finite.  
\item In particular, for the dihedral group $D_4$ with the above presentation we have $\gamma(0) = 1$,  $\gamma(1) = 4$,  $\gamma(2) = 7$ and   $\gamma(n) = 8$ for any $n\geq 3$.
\item The infnite cyclic group $\Z$, generated by $S = \{+1,-1\}$,  has linear grow $\gamma(n) = 2n+1$.  Conversely, any finitely generated group $\Gamma$  whose growth satisfies $\limsup_{n\to \infty}\gamma(n)/n < \infty$ is either finite or contains an infinite cyclic group of finite index. This was proved by   J. Justin in 1971,   see \cite[Theorem 3.1]{Mann} or \cite[Theorem 7.26]{Ceccherini}.
\item  A finitely generated abelian  group of rank $d$ has polynomial growth of degree $d$ for any system of generators, that is $c_1 n^d \leq \gamma(n) \leq c_2 n^d$ for some constants $c_1,c_2$. The precise growth function of  the abelian group $\Z^d$, with respect to its standard set of generators, is given in \cite[page 9]{Harpe2022}. 
\item   The free group $F_q$ on $q$ generators has exponential growth if $q\geq 2$. More precisely, we have:
$$
 \gamma(n) = {\dfrac {q \left( 2\,q-1 \right) ^{n}-1}{q-1}}$$
This  is easily proved by induction from the fact that $\gamma(0) = 1$ and $\gamma(n)-\gamma(n-1) = 2q(2q-1)^{n-1}$ is the number of elements in  $F_q$ at distance exactly $n$ from $e$.
\item   The Heisenberg group $\mathbb{H}_3 \subset GL_3(\Z)$ of $3\times 3$ upper-triangular matrices with $1$ on the diagonal has polynomial growth of degree $4$. This follows from \cite[Theorem 2]{Bass}, which in fact computes the polynomial  growth degree of an arbitrary finitely generated nilpotent group.
\item A celebrated Theorem, proved in 1981 by M. Gromov, states that a finitely generated group has at most polynomial growth if and only if it contains a nilpotent subgroup of finite index. A detailed, self-contained  reference for that result is Chapter 12 in  \cite{Ceccherini}.
\end{enumerate}

\medskip

In general, the growth $\gamma_{_S} (n)$ of an infinite, finitely generated group, is at least a linear and at most an exponential function of $n$.  An important development of Geometric Group Theory has been the discovery in 1983 by  R. Grigorchuk of a class of groups with ``intermediate growth'', that is  finitely generated groups whose growth $\gamma_s$ satisfies
$$
 \lim_{n\to \infty} \;  \frac{\log (\gamma_{_S} (n))}{n} =  0
 \quad \text{ and } \quad 
  \liminf_{n\to \infty} \frac{\gamma_{_S} (n)}{n^d} = \infty,
  \quad  \forall d\in \N.
$$
For an updated presentation on the subject of intermediate-growth groups, we refer to \cite{Grigorchuk} and \cite{Mann}, and the recent papers \cite{Erschler2020b} and \cite{Brieussel}. 

\medskip 
 
\textbf{(iii) The inverse growth function.}  The  \textit{inverse growth function}  \index{Inverse growth} is the function   $\phi_S: \R_+  \to \N \cup \{\infty\}$ defined as  
\begin{equation}\label{def.phi}
 \phi_S(t) = \min \{n \in \N \cup \{0\} \mid \gamma_{_S} (n) \geq t\}.
\end{equation}
This is the smallest function such that $\gamma_{_S} (\phi_S(t)) \geq t$ for any $t\in \R_+ $, equivalently $\phi_S(t)$ is the smallest integer $n$ such that $B_S(n)$ contains at least $t$ elements. Note that $\phi_S(t) = \infty$ if and only if $t > |\Gamma|$, in particular $\phi_s(t)$ is always finite if $\Gamma$ is an infinite group. Observe also that  $\phi_S$ is a left inverse of $\gamma_S$, that is 
$
  \phi_S \left( \gamma_S (n)\right) = n
$
for any integer $n\in \N$.

\medskip 

\textbf{(iv) Various definitions of the boundary of a  subset in $\Gamma$.} 
There are several ways to define a notion of boundary of a  non empty  subset $D \subset \Gamma$  with respect to a symmetric generating set $S\subset \Gamma$.
We first define the \textit{inner boundary}  as
$$
   \partial_S D = \{x\in D \mid \operatorname{dist}_S(x, \Gamma \setminus D) = 1\} 
  = \{x\in D \mid \exists \,  s \in S \  \text{ such that }  sx \not\in D\}.
$$
We also define the \textit{outer boundary} of  $D \subset \Gamma$   as
$$
 \partial_S' D = \{x\in \Gamma \setminus D \mid \operatorname{dist}_S(x, D) = 1\}
=  \{x \in \Gamma \setminus  D \mid \exists s \in S \  \text{ such that }  sx \in D\}.
$$
A third boundary that is sometimes considered is the \textit{edge boundary}, which is the set of edges in the Cayley graph joining a point in $D$ to a point in $\Gamma \setminus D$. We denote it by
$$
  E_S(D) = \big\{ \{x,y\} \mid x\in D, \ y\in \Gamma \setminus D, \ d_S(x,y) =1 \big\}.
$$
Observe that  $(\partial_SD \cup \partial_S'D, E_S(D))$ is a  bipartite 
subgraph of $X(\Gamma,S)$ such that each vertex in $\partial_S D$ is related to at least one and at most $|S|$ vertices in $\partial_S D'$, and vice versa. In particular we have 
$$
 \max\{|\partial_S D|, |\partial'_S D|\}
 \leq  |E_S(D)| \leq |S| \cdot  \min\{|\partial_S D|, |\partial'_S D|\}.
$$

\medskip 

\textbf{(v) The isoperimetric profile.}
The \emph{isoperimetric profile} of $(\Gamma,S)$, introduced by M. Gromov   in \cite{Gromov93},  is the smallest function 
$\mathrm{I}_S : \N \to \N$  such that  $|\partial_S D| \geq \mathrm{I}_s(|D|)$ for any subset $D \subset \Gamma$. More explicitly, it is the function  \index{Isoperimetric profile}
$$
   \mathrm{I}_S(m) = \min \{ |\partial_S D|  \mid  D \subset \Gamma, \ |D| = m \}.
$$

\medskip 

\textbf{(vi)  F{\o}lner sequences and the F{\o}lner function.}  
A \textit{F{\o}lner sequence} \index{F{\o}lner sequence} in the group $\Gamma$ is a sequence of finite nonempty sets  
$D_j \subset \Gamma$ such that 
$$
 \lim_{j\to \infty} \frac{|D_i \setminus xD_i|}{|D_i|} = 0
$$
for any $x\in \Gamma$.  Equivalently, if $S\subset \Gamma$ is a finite generating set, then $\{D_j\}$ is a F{\o}lner sequence if and only if
$$
 \lim_{j\to \infty} \frac{|\partial_SD_i|}{|D_i|} = 0.
$$
It has been proved by E.  F{\o}lner in 1955 that a finitely generated group is amenable if and only if there exists a F{\o}lner sequence in that group,
see e.g. \cite[Corollary 14.24]{Ceccherini}. 

\smallskip

A related notion is the  \textit{F{\o}lner function} \index{F{\o}lner function} of $(\Gamma, S)$, introduced by A. Vershik in \cite{Vershik}. It is the function $\text{F{\o}l} : \N \to \N \cup \{\infty\}$   defined as 
$$
 \text{F{\o}l}(n) = \min\big\{ k \in \N \mid \exists D \subset \Gamma  \text{ s.t. }   |D| = k \text{ and }  k \geq \ n|\partial_SD|\big\},
$$
with the convention that $\min (\emptyset) = +\infty$. In other words, $\text{F{\o}l}(n)$ is the cardinality of the smallest set $D\subset \Gamma$ such that 
$$\frac{|\partial_SD|}{|D|} \leq \frac{1}{n}.$$
Observe also the following relation between the isoperimetric profile and the F{\o}lner function:
$$
    \mathrm{I}_S(\text{F{\o}l}(n)) \leq \frac{\text{F{\o}l}(n)}{n},
$$
for all $n\in \N$ such that $\text{F{\o}l}(n) <\infty$.

\medskip 

We also mention that the F{\o}lner function is sometimes defined in a slightly different way:
$$
 \Phi(n) = \min\big\{ k \in \N \mid \exists D \subset \Gamma  \text{ s.t. }   |D| = k \text{ and } |s^{-1}D \triangle D|\cdot n \leq k, \ \forall s\in S\big\},
$$
where 
$$
 s^{-1}D \triangle D = (s^{-1}D\cup D) \setminus (s^{-1}D\cap D) =
 \{x\in D \mid sx \not\in D\} \cup  \{x\in s^{-1}D \mid x \not\in D\}
$$
 is the symmetric difference. Observe that 
$$
 \frac{1}{|S|} |\partial_SD| \leq  \max_{s\in S} |sD \triangle D| \leq 2 |\partial_SD|,
$$

therefore 
$$
 \Phi\left(\tfrac{n}{2} \right) \leq \text{F{\o}l}(n) \leq \Phi(|S|n).
$$

Some simple examples of F{\o}lner functions are:
\begin{enumerate}[$\circ$]
\item It always holds that $\text{F{\o}l}(1) = 1$ and $\text{F{\o}l}(n) \geq n\, $ if $n \leq  |\Gamma|$. 
\item For a finite group, we have  $\text{F{\o}l}(n) = |\Gamma|$ if $n \geq |\Gamma|$ (because $\partial_S \Gamma = \emptyset$).
\item Returning to the  Dihedral group $D_4$, we have  $\text{F{\o}l}(1) = 1$, $\text{F{\o}l}(2)=7$ and  $\text{F{\o}l}(n)=8$
  for $n \geq 3$. 
\item The F{\o}lner function of $\Z$ is $\text{F{\o}l}(1) = 1$ and  $\text{F{\o}l}(n) = 2n$ if $n \geq 2$.
\item  For the free group on $q\geq 2$ generators $F_q$, we have  $\text{F{\o}l}(n) = \infty$ if $n \geq 2$, this follows from Proposition \ref{prop.caslibre} below. 
\item The group  $\Gamma$ is  amenable if and only if  $\text{F{\o}l}(n) < \infty$ for all integers $n$.
\end{enumerate}

\medskip

The sequence of balls $\{B_S(n)\}$ in an amenable finitely generated group $\Gamma$ is not always a  F{\o}lner sequences, but the following facts are known:

\begin{enumerate}[ \quad  (a)]
\item The sequence of balls $\{B_S(n)\} \subset \Gamma$ is F{\o}lner  for any group of polynomial growth, see  \cite[prop. 12.7]{Mann}.
\item Conversely, if the sequence of balls   $\{B_S(n)\}$ is  F{\o}lner, then $\Gamma$ has subexponential growth,  see  \cite[prop. 12.7]{Mann}.
\item If $|\partial B_S(n)| \geq c |B_S(n)|$ for any $n\in \N$ and some $c>0$, then $\Gamma$ has exponential growth. In fact we clearly have $0 < c < 1$ and $\gamma(n) - \gamma(n-1) \geq c \gamma(n)$. Therefore $\gamma(n)(1-c) \geq \gamma(n-1)$ for any integer $n\geq 1$ and since $\gamma (0) = 1$ we have $\gamma(n) \geq \left(\frac{1}{1-c}\right)^n$.
\item If  $\Gamma$ has subexponential growth, then there exists a sequence $\{n_j\} \subset \N$ such that  $\{B_{n_j}\}$ is F{\o}lner,
 see  \cite[\S 3.4]{Harpe}.
\end{enumerate}

\medskip 

\textbf{(vii)  The Cheeger constant}. 
We finally define the  \emph{Cheeger constant} \index{Cheeger constant}  for a finite group $\Gamma$ with respect to $S$ as follows:
$$
  h_S(\Gamma) = \min \left\{\frac{|E_S(D)|}{|D|} \ \big| \ D\subset \Gamma, \ 0 < |D| \leq \frac{1}{2} |\Gamma| \right\}. 
$$
The Buser--Cheeger inequality for a finite group tells us that 
$$
 \frac{1}{2}\lambda_1  \leq  h_S(\Gamma) \leq  \sqrt{2|S|\lambda_1},
$$
where $\lambda_1$ is the smallest eigenvalue of the combinatorial Laplace operator $L : \R^{\Gamma} \to \R^{\Gamma}$, which is defined as
$$
  Lf(x) = \sum_{d_S(x,y)=1} (f(x)-f(y)).
$$
See Propositions  4.2.4  and 4.2.5 in  \cite{Lubotzky}.

\section{The isoperimetric inequality in finitely generated groups and some consequences}

\subsection{Isoperimetric inequalities in free groups and free abelian groups}

The general isoperimetric inequality for general finitely generated  groups is discussed in the next  subsection. Here we first  consider the special case of  free groups and free abelian groups. The following elementary result  is probably well known to the experts, but we did not find it in the literature:

\begin{proposition}\label{prop.caslibre}
Let $F_q$ be the free group on $q$ generators $\{s_1, \dots, s_q\} \subset F_q$, and consider the symmetric generating set $S=  \{s_1, \dots, s_q,s_1^{-1}, \dots, s_q^{-1}\}$.
Then for any  finite, non empty subset $D \subset F_q$, we have  
\begin{equation}\label{isopFree1}
  \frac{ |\partial D| }{|D|} \geq   \frac{q-1}{q} + \frac{m}{q|D|},
\end{equation}
where $m$ is the number of connected components of $D$. 
Furthermore, if the set $D$ is connected, then the outer boundary satisfies the following identity:
\begin{equation}\label{isopFree2}
  |\partial'D| = (2q-2)|D| + 2.
\end{equation}
\end{proposition}

\smallskip 

By definition, two vertices in a subset  $D\subset F_q$ belong to the same connected component if and only if they  can be joined by an edge path contained in $D$.

\medskip

\textbf{Proof.}  Recall that the  Cayley graph of $F_q$ is the infinite  regular tree of degree $2q$. We first prove the second statement by induction on the  cardinality  of $D$.
If $|D| = 1$, then $D$ contains exactly one point and the outer boundary of $D$ is the set of all neighbors of that point. There are $2q$ such neighbors,  parameterized by $S$,  therefore Equation \eqref{isopFree2} is  trivially satisfied.  Assume now that  $D$  is connected and contains at least two points  and chose a base point $x_0$ 
in $X$ and a point $x\in D$ at maximum distance from $x_0$  (such a point is sometimes called a \textit{leaf} in $D$). Then $x$ has one neighbor in $D$ and
$2q-1$ neighbors outside $D$, denote them by $\{y_1, \dots ,y_{2q-1}\}$.  Clearly, the outer boundary of $D^- = D \setminus \{x\}$ is given by 
$$
 \partial' D^- = \left(\partial' D \setminus \{y_1, \dots ,y_{2q-1}\}\right)  \cup \{x\}.
$$
Therefore  we have
$$
  |D| = |D^-|+1 \quad \text{and} \quad |\partial' D| =   |\partial'D^-|  + (2q-2),
$$ 
and we conclude the proof of \eqref{isopFree2}  by induction. 

\smallskip

The proof of the first statement follows immediately for a connected set $D\subset F_q$ from  \eqref{isopFree2}  and the fact that $ |\partial D|  \geq \frac{1}{2q}  |\partial' D| $.
If $D$ has several connected components $D_1, \dots, D_m$, then  the inner boundary $\partial D$ is the disjoint union  of the  $\partial D_j = D_j \cap \partial D$ for any $j = 1, \dots, m$ and we have 
$$
 (2q-2)|D| + 2m = \sum_{j=1}^m \left( (2q-2)|D_j| + 2 \right)  =  \sum_{j=1}^m |\partial' D_j| \leq 
 2q \sum_{j=1}^m   |\partial D_j| = 2q  |\partial D|,
$$
which implies the  inequality \eqref{isopFree1}.
 
\qed

\medskip 

Note that the reason why \eqref{isopFree2} may fail for a disconnected set is that the  outer boundaries of two different connected components of that set may not be disjoint.
Regarding  free abelian groups, we have:
\begin{proposition}\label{prop.caslibreabelien}
For any  finite, non empty,  subset $D \subset\Z^d$ we have  
\begin{equation}\label{isopFreeab}
 \frac{ |\partial D| }{|D|} \geq  \frac{1}{|D|^{1/d}} .
\end{equation}
\end{proposition}

\textbf{Proof.} 
The results follows from the  \textit{Loomis-Whitney Inequality} \cite{LW},  which states that for any  finite subset $D\subset \Z^d$ 
we have 
$$
  |D|^{d-1} \leq    \prod_{j=1}^d |\pi_j(D)|,
$$
where   $\pi_j : \Z^d \to \Z^{d-1}$ is the projection in direction of the $j^{\text{th}}$ coordinate. 
We thus  have
$$
  |\partial D| \geq \max_{1\leq j \leq d}{|\pi_j(D)|} \geq \left(\prod_{j=1}^d |\pi_j(D)|\right)^{\frac{1}{d}} \geq |D|^{\frac{d-1}{d}}.
$$
\qed

\medskip

Inequality  \eqref{isopFreeab} can be improved for the outer boundary, see \cite[Theorem 6.22]{LyonsPeres}.

\subsection{Statement of the Main Result}

The heart of the argument in Gromov's proof of Inequality \eqref{ineqG} involves an estimate of the number of elements in a finite set $D \subset \Gamma$ that are moved outside $D$ by the action of a ball $B(r) \subset \Gamma$ whose radius $r$ is well chosen. To quantitatively express this estimate in terms of the growth function of the group, it is convenient to introduce the following notion: 

\begin{definition} \rm  
Given an arbitrary subset $E \subset \R_+ = [0, +\infty)$ and a  function $g : E \to \R_+$, we define a new function $\U_{E,g} : \R_+ \to \mathbb{R}\cup \{\infty\}$ by
\begin{equation}\label{defUT}
   \U_{E,g}(t)  = \sup \left\{\frac{1}{r}  \left(1- \frac{t}{g(r)} \right) \mid  r\in E \setminus \{0\} \right\}.
\end{equation}
\index{$\U$-transform} We will call this function the $\U$-\textit{transform} of $g$. When the domain $E \subset \R_+$ is fixed, we  usually  write $\U_{g}(t)$ instead of $\U_{E,g}(t)$.
\end{definition}

\medskip 

 As a first example, we mention that  the $\U$-transform of the polynomial function $g(r) = (d+1)r^d$ is the function $\U_{g}(t) = \frac{d}{d+1}t^{-1/d}$.   More examples are given in  Section \ref{onUtransform}, where some
computational techniques  and   basic properties of the  $\U$-transform will be given.
 We will in particular explain in \S \ref{sec.legendre} that the $\U$-transform  is nothing else than  a variant of the classical Legendre transform obtained by some change of variables.

\medskip 

Using the $\U$-transform, we now state the main result of the present paper, which is  the following version of  isoperimetric inequality in finitely generated groups:

\begin{theorem} \label{maintheoremA}
In an arbitrary  group $\Gamma$, generated by the finite symmetric set $S = S^{-1} \subset \Gamma$,
the following isoperimetric  inequality  holds for any non empty finite subset $D \subset \Gamma$:
\begin{equation}\label{MainIsopIneq}
  \frac{ |\partial_SD|}{ |D|} \geq  \U_{\gamma_{_S}} (|D|),
\end{equation}
where the growth function $\gamma_{_S} : \N \to \N$ and the boundary $\partial_{S} D$ are defined with respect to the generating set $S$. 
Furthermore, the following holds for the F{\o}lner function of $(\Gamma,S)$:
\begin{equation}\label{FMainIsopIneq}
 \U_{\gamma_{_S}}(\textrm{\rm F{\o}l}(n)) \leq \frac{1}{n}.
\end{equation}
If $\Gamma$ is a finite group, then its Cheeger constant satisfies
\begin{equation}\label{CheegerIsopIneq}
  h_s(\Gamma) \geq  \U_{\gamma_{_S}}\left(\frac{1}{2} |\Gamma|\right).
\end{equation}
\end{theorem}

\medskip 

We will prove this Theorem  in the next subsection. We first derive the following result, which is stated in \cite{PittetStankov} and  \cite{Santos},  and  implies the inequality of Coulhon and Saloff-Coste:
 
\medskip

\begin{corollary} \label{maincorollary}
Let $\Gamma$ be a finitely generated group. For any non empty finite subset $D\subset \Gamma$ and any 
real number $\lambda > 0$, we have
\begin{equation}\label{mainineq}
   \frac{ |\partial_SD|}{ |D|}   \geq  \left(1 - \frac{1}{\lambda}\right) \frac{1}{\phi_S(\lambda |D|)}.
\end{equation}
\end{corollary}

Choosing $\lambda = 2$ in this inequality gives us \eqref{ineqG}.

\medskip

\textbf{Proof.}  
If $\lambda |D| > |\Gamma|$, then $\phi_S(\lambda |D|) = \infty$ and there is nothing to prove. We thus assume  that 
 $\lambda |D| \leq  |\Gamma|$ and set  $r = \phi_S(\lambda |D|)$. Then  $\lambda |D|  \leq  \gamma_{_S} (r)$ and  Theorem \ref{maintheoremA} implies that 
$$
   \frac{ |\partial_SD|}{|D|}    \geq \U_{\gamma_{_S} }(|D|) \geq   \frac{1}{r} \left( 1- \frac{|D|}{\gamma_{_S} (r)}\right) 
 \geq \left( 1 - \frac{1}{\lambda}\right) \frac{1}{\phi_S(\lambda |D|)}.
$$
\qed

\subsection{Proof of the Main Theorem} \label{secproofm}

Theorem \ref{maintheoremA} will be a direct consequence of  the following  stronger result:

\begin{theorem}\label{maintheorem}
Given an arbitrary  group $\Gamma$ generated by the finite symmetric set $S = S^{-1} \subset \Gamma$,
the following isoperimetric inequality  holds for any finite, non empty  subset $D \subset \Gamma$:
\begin{equation}\label{MainIsopIneqStrong}
  \frac{ |\partial_SD|}{ |D|}   \geq   
  \sup_{r \in \N} \left( \frac{\gamma_S(r) - |D|}{r\gamma_S(r)  - \sum_{k=0}^{r-1}\gamma_S(k)}\right),
\end{equation}
where the growth function $\gamma_{_S} : \N \to \N$ and the boundary $\partial_{S} D$ are defined with respect to the
generating set $S$. 
\end{theorem}

\medskip 

We will use the following notation: 
For any $k\in \N$, we   denote by $S_S(k) = \{x \in \Gamma \mid \|x\|_S =k\}$ the sphere of radius $k$ in $\Gamma$,  
and we set $\sigma_S(k) = |S_S(k)|$. Note that the  ball $B_S(r)$ is the disjoint union of the spheres $S_S(k)$ for 
$0 \leq k \leq n$. In particular we have 
\begin{equation}\label{eqgammasigma}
  \gamma_S (r) = \sum_{k=0}^r\sigma_S(k)\quad \text{ and }  \quad \sigma_S(k) = \left( \gamma_S (k) - \gamma_S (k-1)\right).
\end{equation}

\medskip 

\textbf{Proof of Theorems \ref{maintheorem}.}
The proof follows a strategy similar to that in Gromov's book \cite[\S 33]{Gromov}, see also \cite[\S 6.7]{LyonsPeres}.
We first claim that for any $y\in S(k)$, we have
\begin{equation} \label{eqb1}
   |\{x\in D \mid yx \not\in D\}|  \leq k\,|\partial D|. 
\end{equation}
Indeed, if $k=0$ then $y = \id_{\Gamma}$ and the claim is trivial. Moreover, the following inclusion holds for any $y\in \Gamma$ and any $s\in S$ :
\begin{equation} \label{inclusionbdy}
 \{x\in D \mid syx \not\in D\} \subset \{x\in D \mid yx \not\in D\}  \cup \{x\in D \mid yx \in\partial_SD\},
\end{equation}
thus \eqref{eqb1} follows by induction on $k$. Let us now   set 
$
  P_r(D) =  \{ (x,y ) \in D\times B_S(r)  \mid  yx \notin D  \},
$
we then have from  \eqref{eqb1}:
\begin{equation}\label{eqb2}
  |P_r(D) |  =  \sum_{k=0}^r \ | \{ (x,y ) \in D\times S_S(k)  \mid  yx \notin D  \}|   \leq  |\partial D| \, \sum_{k=0}^r k \, \sigma(k).
\end{equation}
On the other hand,  the following  inequality is obvious for every $x\in D$ and  any $r\in \mathbb{N}$: 
$$
  |(\{y\in B_S(r) \mid  yx \notin D \} |   =  |\{y \in B_S(r) \mid  y  \notin Dx^{-1} \}|   \geq  |B_S(r)| - |D|.
$$
This  can be written as
\begin{equation}\label{eqb3}
  \gamma_S(r) - |D|  \leq  | \{y\in B_S(r) \mid  yx \notin D \}|.
\end{equation}
From  \eqref{eqb2} and  \eqref{eqb3}, one obtains the following inequalities  for any $r\in \N$:
$$
  |D| \left(\gamma_S(r) - |D|\right)  \leq  |P_r(D)| \leq  |\partial D| \, \sum_{k=0}^r k \, \sigma_S(k),
$$
from which the inequality 
\begin{equation}\label{eq.rass}
  \frac{ |\partial_SD|}{ |D|}   \geq   \sup_{r \in \N} \left( \frac{\gamma_S(r) - |D|}{\sum_{k=1}^r k \sigma(k)}\right)
\end{equation}
follows immediately. Inequality \eqref{MainIsopIneqStrong} follows  now from  \eqref{eq.rass} and the obvious identity
$$
  \sum_{k=1}^r k \sigma(k) =  r\gamma_S(r)  - \sum_{k=0}^{r-1}\gamma_S(k).
$$
\qed

\medskip

\textbf{Proof of Theorem \ref{maintheoremA}.}
The inequality \eqref{MainIsopIneq}   follows now immediately from \eqref{MainIsopIneqStrong} and the definition \eqref{defUT} of the $\U$-transform:
$$
  \frac{ |\partial_SD|}{ |D|}   \geq   
  \sup_{r \in \N} \left( \frac{\gamma_S(r) - |D|}{r\gamma_S(r)  - \sum_{k=0}^{r-1}\gamma_S(k)}\right)
  \geq   \sup_{r \in \N} \left( \frac{\gamma_S(r) - |D|}{r\gamma_S(r) }\right) = \U_S|D|).
$$ 
Inequalities \eqref{FMainIsopIneq} and  \eqref{CheegerIsopIneq} are immediate consequences of \eqref{MainIsopIneq} and the definitions of the F{\o}lner function and the Cheeger constant. 

\qed

\subsection{Some Consequences of the Main Result} \label{sec.consequences} 

In this section and the next one, we derive some consequences of  Theorem \ref{maintheoremA}. The proofs  use the basic properties of the $\U$-transform developed in Section \ref{onUtransform}. We begin with the following  statement on groups with polynomial growth:
\begin{corollary} \label{polynomialgrowth}
Let $\Gamma$ be a finitely generated group whose  growth function $\gamma_{_S} $ satisfies  $\gamma_{_S} (n-1) \geq   Cn^d$ for some constants $C >0$ and  $d \geq 1$ and any integer $n \geq 1$. Then the following isoperimetric inequality 
\begin{equation}\label{isoppol}
    \frac{ |\partial_SD|}{|D|}  \geq  \frac{C^{\frac{1}{d}} d}{(d+1)^{1+\frac{1}{d}}}|D|^{-\frac{1}{d}}
\end{equation}
holds for any finite, non empty subset $D\subset \Gamma$.

\end{corollary}

\medskip

\textbf{Proof.} The inequality \eqref{isoppol} follows from Theorem \ref{maintheoremA}   together with Lemma \ref{URZ} and the computation \eqref{Uofrd}  of the polynomial $\U$-transform given in Example \ref{premierexample} below.

\qed

\medskip

By comparison, for the same growth function, the inequality \eqref{ineqG} gives us the estimate:
$$
    \frac{ |\partial_SD|}{|D|}  \geq  
\frac{C^{\frac{1}{d}}}{2^{1+\frac{1}{d}}\left(1-\left(\frac{C}{2|D|}\right)^{\frac{1}{d}}\right)}|D|^{-\frac{1}{d}}.
$$
Note that for  large $d$, the constant in the latter inequality is about one half that in \eqref{isoppol}.

For group with exponential growth, we have the following 
\begin{corollary} \label{exponenialgrowth}
Let $\Gamma$ be a finitely generated group whose  growth function $\gamma_{_S} $ satisfies  $\gamma_{_S} (n-1) \geq   C \ee^{bn^{\alpha}}$ for some constants  $0 < \alpha \leq 1$, \  $C, b>0$,  and any integer $n \geq 1$. Then we have
\begin{equation}\label{isoexp}
    \frac{ |\partial_SD|}{|D|}  \geq  \left( \frac{b}{\log(|D|)  + o(\log(|D|))}\right)^{1/{\alpha}}
\end{equation}
for any finite subset $D\subset \Gamma$  with at least two elements.
\end{corollary}

\medskip 
 
\textbf{Proof.} The result follows from Theorem \ref{maintheoremA} combined with the inequality \eqref{ineqUofexppp}   in   Example   \ref{secondexample}  below.

\qed 

\medskip

Note that for the same growth function,  \eqref{ineqG} gives us
$$
    \frac{ |\partial_SD|}{|D|}  \geq  \frac{1}{2^{1/\alpha}}  \left( \frac{b}{\log(|D|)  + \log(2/C)}\right)^{1/{\alpha}}.
$$

\bigskip

Our next result gives   lower bounds for the F{\o}lner function:
\begin{corollary}\label{corfolner}
Let $\Gamma$ be a finitely generated group with  growth function $\gamma_{_S} $.
\begin{enumerate}[(a)]
    \item If $\gamma_{_S} (n-1) \geq C n^d$ for some $d \geq 1$, then \ 
$\displaystyle  {\text{\rm F{\o}l}(n) } \geq \frac{C d^d }{(1+d)^{1+d}}\cdot  n^d$,
    \item If $\gamma_{_S} (n-1) \geq C \ee^{bn}$ for some $b>0$, then \
$\displaystyle  \text{\rm F{\o}l}(n) \geq  \frac{C}{\alpha e b}\cdot \frac{\exp\left({bn^{\alpha}}\right)}{n^{\alpha}} $.
\end{enumerate}
\end{corollary}

\medskip

\textbf{Proof.}  Using \eqref{FMainIsopIneq} and the calculation in Example \ref{premierexample}, we see that if $\gamma_{_S} (n-1) \geq C n^d$, then
$$
 \frac{1}{n} \geq \U_{\gamma_S}(\text{\rm F{\o}l}(n)) =  \frac{d C^{\frac{1}{d}}  }{(d+1)^{1+\frac{1}{d}}}\cdot  \frac{1}{\left( \text{\rm F{\o}l}(n)\right)^{\frac{1}{d}}},
$$
which proves (a). 

\smallskip

To prove (b), we use the calculation in Example \ref{secondexample}. In particular if  $\gamma_{_S} (n-1) \geq  g(n)  = C \ee^{bn^\alpha}$, 
then Equation \eqref{InvUofexpr} with $u = 1/n$ and $t = \text{\rm F{\o}l}(n)$ implies that 
$$
\U_{\gamma_S}(\text{\rm F{\o}l}(n)) \leq  \frac{1}{n}   \quad \Rightarrow \quad
 \text{\rm F{\o}l}(n) \geq  \frac{C}{\alpha e b}\cdot \frac{\exp\left({bn^{\alpha}}\right)}{n^{\alpha}} 
$$
\qed

\subsection{Asymptotic Estimates}  

In this section we formulate some asymptotic estimates on the isopermietric ratio and the F{\o}lner function   for some groups with non polynomial growth. We will need the following 
somewhat technical definition:

\begin{definition} \rm \label{def.TSPG}
We will say that a function $g : \R_+ \to \R_+$ has \textit{Tame Superpolynomial Growth}, abbreviated as (TSPG), if it is everywhere 
differentiable,  with $g'>0$ and 
\begin{equation}\label{TSPG}
 \lim_{r\to \infty} \left(\frac{r g'(r)  g(\lambda r) }{g(r)^2}\right)  = 
 \begin{cases} 0, & \text{if } \ 0 < \lambda < 1, \\ \infty, & \text{if } \  \lambda = 1.  \end{cases}
 \tag{TSPG}
\end{equation}
\end{definition}
This condition can   also be written for $f(r) = \log(g(r))$ as follows:
$$
 \lim_{r\to \infty} r f'(r) = \infty \quad \text{and} \quad   \lim_{r\to \infty} \frac{r f'(r)}{\exp({f(r)-f(\lambda r)})} = 0, 
 \ (\text{for any } 0 < \lambda < 1).
$$
An obvious example   is the exponential function. Other examples are the functions $g_1(r) =  \exp({ar^{\beta}})$ with $a >0$ and $0 < \beta \leq 1$, \ $g_2(r) =  \exp \left(\frac{ar}{\log(r+1)^{\alpha}} \right)$ for some $a>0$ and any $\alpha \in \R$ and $g_3(r) = r^{\sqrt{r}}$.

\smallskip

This notion is interesting because the non polynomial  group growths that are described in the literature are bounded below by functions satisfying (TSGP), see e.g. \cite{Grigorchuk,Mann}. Our next results are formulated under this hypothesis:

\begin{corollary}  \label{cortomain}
Let $\Gamma$ be a finitely generated group whose  growth function $\gamma_{_S} $ satisfies  $\displaystyle \gamma_{_S} (r-1) \geq g(r)$, for any $r \geq 1$, where  $g: \R_+ \to \R_+$ satisfies the growth condition (\ref{TSPG}). 
Then for  any $\varepsilon>0$, there exists $N = N(\varepsilon) \in \N$  such that the following inequality holds
\begin{equation}\label{isopexp}
     \frac{ |\partial_SD|}{|D|}   \geq   \frac{1-\varepsilon}{g^{-1}(|D|)},
\end{equation}
for any finite subset $D\subset \Gamma$ such that $|D| \geq  N(\varepsilon) $.
\end{corollary}

\textbf{Proof.} Inequality \eqref{isopexp} follows from Theorem   \ref{maintheoremA}, together with Lemma \ref{URZ} and
Proposition \ref{prop.subexpg} below. 

\qed
 
\medskip

Using the previous result, we obtain the following asymptotic estimate for the F{\o}lner function of groups of intermediate growth, thus completing Corollary \ref{corfolner}. 
\begin{corollary}\label{Corest.Folner2}
Let $\Gamma$ be an infinite amenable group  satisfying the hypothesis of  Corollary \ref{cortomain}. Then for  any $\varepsilon>0$, there exists $N = N(\varepsilon) \in \N$  such that for any $n \geq N$,
\begin{equation}\label{est.Folner2}
 \text{\rm F{\o}l}(n) \geq  g((1-\varepsilon) n).
\end{equation}
In other words, we have  
\begin{equation}\label{est.Folner3}
  \liminf_{n\to \infty}  \,  \frac{g^{-1}(\textrm{\rm F{\o}l}(n))}{n}\geq 1.
\end{equation}
\end{corollary} 

\textbf{Proof.} Suppose \eqref{est.Folner2} does not hold, then there exists $\eta > 0$ and a sequence ${n_j} \subset \N$ such that $n_j \to \infty$ and 
${g^{-1}(\textrm{F{\o}l}(n_j))}/{n_j}\leq (1 - \eta)$ \, for any  integer $j$. From the definition of the F{\o}lner function, one can then find finite subsets $D_j \subset \Gamma$ such that $|D_j| =\textrm{F{\o}l}(n_j)$ and 
$$
  \frac{|\partial_S D_j|}{|D_j|} \leq \frac{1}{n_j} \leq \frac{1-\eta}{g^{-1}(\text{F{\o}l}(n_j))}  = \frac{1-\eta}{g^{-1}(|D_j|)}.
$$
This inequality  contradicts  Corollary \ref{cortomain} since $\displaystyle \lim_{j\to \infty} |D_j| = \infty$.

\qed

\medskip 

\begin{remark} \rm 
 Using the  other version $\Phi$ of the F{\o}lner function, L. Bartholdi gave a direct proof of the following inequality:
 $$
  \Phi(n) \geq \frac{1}{2} g(n),
 $$
 for any $n\in \N$, see \cite[page 455]{Bartholdi}. Because $\text{F{\o}l}(n) \geq \Phi\left(\frac{n}{2} \right)$, the above inequality implies 
 $$
  \text{F{\o}l}(n) \geq \frac{1}{2}g\left(\frac{n}{2} \right),
 $$
which is slightly better than  Proposition 14.100 in \cite{Ceccherini}.  Note that Corollaries  \ref{corfolner}  and  \ref{Corest.Folner2}  improve this inequality. 
\end{remark}

\medskip

The following question is inspired by the recent work of C. Pittet and B. Stankov in \cite{PittetStankov,Stankov},  and the previous Corollary.
\begin{question} \label{QKg} \rm 
Given a smooth, unbounded monotone increasing function $g: \R\to\R$, we ask for the asymptotically smallest possible F{\o}lner function among all groups $\Gamma$ with a finite, symmetric, generating set $S$ whose growth function satisfies  $\gamma_{_S}  (n-1) \geq g(n)$. Specifically, given the function $g$, we ask for the value of 
\begin{equation}\label{defKg}
    K_g = \inf\left\{\left.\liminf_{n\to \infty}  \,  \frac{g^{-1}(\textrm{\rm F{\o}l}(n))}{n} \  \right| \exists  (\Gamma,S)   \text{ with growth  } \gamma_{_S} (n-1) \geq g(n) \right\}.  
\end{equation}
\end{question}   

For the  polynomial case $g(r) = Cr^d$, Corollary \ref{corfolner} (a) gives us the following lower bound:  
$$
 K_{Cr^d}  \geq  \frac{d}{(1+d)^{1+1/d}}.
$$
On the other hand  Corollary \ref{Corest.Folner2} implies that if for any function  $g$  satisfying  (\ref{TSPG}) we have
$$
  K_g \geq 1.
$$
See also the more explicit  Corollary \ref{corfolner} (b)  for the exponential case.  

\medskip

To find an upper bound for $K_g$, one needs to estimate the size of a  F{\o}lner sequence in an amenable group with growth   $\geq g$.  In \cite{Stankov}, B. Stankov  gave an example where $1 \leq K_g \leq 2$. More precisely, he considers the lamplighter group, that is the  wreath product $\Gamma = \Z \wr (\Z/2\Z)$, with a well chosen generating set $S$ and he computes  in this case that $\gamma_{_S} (n-1) \geq 2^{n-1}$ and 
$\textrm{F{\o}l}(n)\leq 4^n$. It follows that  $g^{-1}(t) = \log_2(2t) = \log_2(t)+1$ and thus
$$
 \frac{g^{-1}(\textrm{F{\o}l}(n))}{n} \leq 2 + \frac{1}{n}.
$$
We refer to \cite{PittetStankov} and \cite{Stankov} for additional  results and related questions. 
 

\subsection{Final  Remarks}

We conclude this Section  with  a few  specific remarks:

\smallskip 

$\circ$ \ The proof of  Theorem \ref{maintheorem}, and therefore all the estimates in the paper, also  holds if we replace the inner  boundary $ \partial_SD$ with the outer boundary  $\partial_S' D$. We only need to change the inclusion \eqref{inclusionbdy} by 
\begin{equation}\label{inclusionbdy2}
 \{x\in D \mid  syx\not\in D\} \subset \{x\in D \mid yx \not\in D\}  \cup \{x\in D \mid syx \in\partial_SD'\},
\end{equation}
and follow the same argument. 

\medskip 

$\circ$ \ The inverse growth function $\phi_S$ defined in \eqref{def.phi} is a variant of the function defined in  \cite{PSC} and \cite{CSC}. These authors use instead the function 
\begin{equation}\label{def.phi2}
 \widetilde{\phi}_S(t) = \min \{n \in \N \mid \gamma_{_S} (n) > t\}.
\end{equation}
Observe that $\widetilde{\phi}_S(t) \geq {\phi}_S(t)$, therefore the isoperimetric inequality \eqref{mainineq} still holds if one uses
$\tilde{\phi}_S$ instead of ${\phi}_S$, and it is in fact a slightly weaker statement in that case.

\medskip 

$\circ$ \ For further results and updated references on the F{\o}lner function and the isoperimetric profile, we refer to Chapter 14 of the book  \cite{Ceccherini}, especially sSctions 14.10 to 14.13. In particular an  upper bound for the F{\o}lner function of a finitely generated nilpotent group in terms of its growth function can be deduced from the proof of  Theorem 14.102 in that book.

\section{Some Properties and Calculations of  the $\mathcal{U}$-transform}  \label{onUtransform}

In this final Section, we present some basic facts on the $\U$-transform and perform some computations that are used in
section \ref{sec.consequences}. As explained  in \S \ref{sec.legendre} below, the $\U$-transform is equivalent to the Legendre transform after some change of variables. However we find it convenient to give direct proofs of the  properties of the $\U$-transform we shall use, rather than obtaining  them as consequences of properties of the Legendre transform.

\smallskip  

Recall   that the $\U-$transform of the function  $g : E \subset  \R_+ \to \R_+$ is defined as
$$
   \U_{E,g}(t)  = \sup \left\{\frac{1}{r}  \left(1- \frac{t}{g(r)} \right) \mid  r\in E \setminus \{0\} \right\}.
$$

\subsection{Basic Properties} 

We first gather some elementary facts:
\begin{lemma}\label{thelemma}
The $\U$-transform  satisfies the following properties:
\begin{enumerate}[(i)]
\item The function $t \mapsto \U_{E,g}(t)$ is non-increasing.
\item If $g : E \to \R_+$ and $E' \subset E$, then $\U_{E',g}(t) \leq \U_{E,g}(t)$ for any $t \in \R_+$.
\item If $g_1, g_2 : E \to \R_+$ and $g_1(r) \leq g_2(r)$ for any $r\in E$,  then $\U_{E,g_1}(t) \leq \U_{E,g_2}(t)$ for any $t \in \R_+$.
\item If $t < \sup_E g(r)$, then $\U_{E,g}(t) > 0$. 
\item If $g$ is increasing and $t > \inf_E g(r)$, then $\U_{E,g}(t) < \infty$. 
\item If  $g,h : \R_+ \to \R_+$ are  any functions such that $h(r) = cg(br)$ for some positive constants $b,c$, then 
$$
  \U_{\R,h}(t) = b \cdot \U_{\R,g}(t/c).
$$
\end{enumerate}
\end{lemma}

\medskip 

\textbf{Proof.}  The first four statements are easy consequences of the definition. Statement (v) follows from the fact that  if $\U_{E,g}(t) = \infty$, then there exists a sequence $\{r_i\} \subset E\setminus \{0\}$ such that $r_i \to 0$  and 
$g(r_i) > t$, contradicting the hypothesis $\lim g(r_i) = \inf g(r_i) < t$. 
To prove the last statement, set $s = br$ and start with  the identity \  
$$
 \frac{1}{r}\left(1- \frac{t}{h(r)}\right) =  \frac{1}{r}\left(1- \frac{t}{cg(br)}\right) =  \frac{b}{s}\left(1- \frac{t/c}{g(s)}\right).
$$
Taking the sup over $r$ on the left hand side and over $s$ on the right hand side proves (vi).

\qed

\medskip

The following Lemma is useful when comparing the $\U$-transform of a function defined on the natural numbers to that of a function
defined on the whole set of positive real numbers.

\begin{lemma}\label{URZ}
Suppose that $g : \R_+ \to \R_+$ and $\gamma : \N  \to \N$ are two functions such that $g$ is 
non decreasing  and   $\gamma(k-1) \geq g(k)$ for any  $k\in \N$, then  \  $\U_{\N,\gamma}(t) \geq \U_{\R,g}(t)$  \ for any $t \in \R_+$.
\end{lemma}

\textbf{Proof.} 
The statement  follows from the fact that  the integer part $k = \lfloor {r} \rfloor$ of any real number  $r$
satisfies  $k \leq r <  k+1$. From  our hypothesis we then have $g(r) \leq g(k+1) \leq \gamma(k)$. Therefore we have
\begin{equation*}\label{ineqRZ}
  \frac{1}{r}  \left(1- \frac{t}{g(r)}\right)  \leq  \frac{1}{k}  \left(1- \frac{t}{\gamma(k)}\right),
\end{equation*}
and the lemma follows immediately.

\qed

\subsection{Computing the $\U$-transform}    

To compute the $\U$-transform of a differentiable function $g : [0, \infty) \to \R_+$, one fixes  $t >  g(0)$ and sets   $f(r) = \frac{1}{r}\left(1- \frac{t}{g(r)}\right)$.
We observe that  $f$ achieves its maximum on $\R_+$ since $f(r) \leq  0$ if $g(r) \leq t$, $f(r) > 0$ if $g(r) > t$ and $\displaystyle \lim_{r\to \infty} f(r) = 0$.
Define now the   function $\rho : (g(0), \infty) \to \R$ by 
\begin{equation}\label{defrho}
 \rho(t) = \max \{r \in \R_+ \mid f(r) = \U_g(t) \},
\end{equation}
and observe that 
\begin{equation}\label{solveU}
  \U_g(t) =  \frac{1}{\rho(t)}  \left(1- \frac{t}{g(\rho(t))}\right).
\end{equation}
We  have thus   reduced the computation of  the $\U$-transform of $g$ to that of $\rho(t)$.  
 
\begin{lemma}\label{taurho}
The function $\tau : \R_+ \to \R_+$ defined by 
\begin{equation}\label{def.tau}
 \tau(r) = \frac{g(r)}{1+r \frac{g'(r)}{g(r)}}
\end{equation}
is a left inverse of $\rho$. 
\end{lemma}

\textbf{Proof.}  
At any point $r$ where  the function $f$ achieves its maximum we have 
\begin{equation*}\label{deriveef}
  f'(r) =  \dfrac{t}{r^2 g(r)}\left(1 + \dfrac{r g'(r)}{g(r)}   - \dfrac{g(r)}{t} \right) = 0.
\end{equation*}
This condition can be written as  $g(r) = t\left( {1+r \frac{g'(r)}{g(r)}}\right)$, or equivalently as $t = \tau(r)$.
In particular we have $\tau (\rho (t)) = t$.

\qed 

\bigskip

\noindent\fbox{%
    \parbox{15.5cm}{%
The procedure to compute the $\U-$transform of a differentiable function $g$ is then as follows:
\begin{enumerate}[(i)]
    \item Compute the function $\tau(r)$ from \eqref{def.tau}.
    \item Compute   $r = \rho(t)$ by solving $\tau(r) = t$.
    \item $\U_g(t)$ is then given by \eqref{solveU}.
\end{enumerate}
}}%

\medskip 

The procedure works fine  provided the function $\tau(r)$ is injective, which will be the case in all the examples we consider.

\medskip 

\begin{example}\label{premierexample} \rm 
As a first example, we consider the  function $g(r)= c r^d$, then $\tau(r) = \dfrac{c}{d+1} r^d$, therefore $\rho(t) = \left(\frac {\left( d+1 \right) }{c} t\right)^{1/d}$ and we have 
\begin{equation}\label{Uofrd}
 \U_g(t) =  \frac{1}{\rho(t)}\left(1-\frac{t}{g(\rho(t))}\right) = \left(\frac{d c^{\frac{1}{d}}  }{(d+1)^{1+\frac{1}{d}}}\right)\cdot t^{-\frac{1}{d}}.
\end{equation}
\end{example}

\bigskip 

Our next example will be to compute the $\U$-transform of exponential functions. The computation is tricky and it will be convenient to introduce the auxiliary function 
$f: [e,\infty) \to [1,\infty)$ defined as
\begin{equation}\label{eq.funf} 
f(x) = \log(x) + \log(f(x))
\end{equation}
The function $f$ is well-defined; it is the  inverse function  of $y \mapsto e^y/y$ (for $y \geq 1$ and $x \geq e$). Moreover it is monotone increasing with  $f(e) = 1$  and $\lim_{x \to \infty} f(x) = \infty$.
It is clear that $f(x) \geq \log(x)$, but more precisely we claim that
\begin{equation}
 f(x) = \log(x) + \log(\log(x)) + o(1).
\end{equation}
To see this, let us write 
$$f(x) = a(x) \log(x),$$
then
$$
  a(x) = \frac{f(x)}{\log(x)} =  \frac{f(x)}{f(x) - \log(f(x))}.
$$

This function is positive, decreasing, and $\displaystyle \lim_{x\to \infty} a(x) = 1$. Therefore  $a(x) = 1+o(1)$ and we have 
\begin{eqnarray*}
f(x)&=& \log(x) + \log(f(x)) = \log(x) + \log \left( a(x) \log(x) \right)
\\ & = & \log(x) + \log \left( \log(x)\right) + o(1).
\end{eqnarray*}
as claimed above. 

\medskip

\textbf{Remark.}  We can also express $f$ in terms of the \emph{Lambert $W$-function}:  \index{Lambert $W$-function}
$$f(x) = -W(-1/x).$$
The function $W$  is defined to be the  inverse  of the function $h(y) = y\ee^y$. Since $h$ is not injective, $W$ is multivaluated,  and here we appeal to  the branch (sometimes denoted as $W_{-1}$) corresponding to the ranges $-1/\ee<x<0$  and $-\infty < y < -1$.

\medskip 

\begin{example}\label{secondexample} \rm 
We now  compute  the $\U$-transform of the   function   $g(r) = ce^{br^{\alpha}}$. 
We   have  $\tau(r) = \frac{ce^{br^{\alpha}}}{1+\alpha b r^{\alpha }}$, and the relation $\tau(\rho(t)) = t$ is thus  equivalent to
$$
 cb \rho(t)^{\alpha}  = \log(t) + \log\left(1+ \alpha b \rho(t)^{\alpha} \right).
$$
This identity can be written as
$$
  \left(b \rho(t)^{\alpha} + \frac{1}{\alpha} \right)  =   \log(\lambda t)  +   \log   \left(b \rho(t)^{\alpha} + \frac{1}{\alpha} \right) ,
$$
with $\lambda = \frac{\alpha}{c}  \ee^{1/\alpha}$.  The solution is given by 
$$
 b \rho(t)^{\alpha} + \frac{1}{\alpha}  = f(\lambda t),
$$
where $f(x)$ is the function defined in \eqref{eq.funf}.
From \eqref{solveU}, one then obtains 
$$
 \U_g(t) = \frac{1}{\rho(t)} \left(1-\frac{t}{ce^{b\rho(t)^{\alpha}}}\right) = \frac{\alpha b\rho(t)^{\alpha-1}}{1+\alpha b\rho(t)^{\alpha}}
=
\frac{ b^{1/\alpha} \left(f(\lambda t)-\frac{1}{\alpha} \right) }{f(\lambda t)},
$$
which be written as 
\begin{equation}\label{Uofexppp}
 \U_g(t) =  \left( \frac{b}{f(\lambda t)}\right)^{\frac{1}{\alpha}} \left(1 -\frac{1}{\alpha f(\lambda t)} \right)^{1-\frac{1}{\alpha}}.
\end{equation}
 
Note that if $0 < \alpha \leq 1$, then $\left(1 -\frac{1}{\alpha f(\lambda t)} \right)^{1-\frac{1}{\alpha}} > 1$ and we have therefore
\begin{equation}\label{ineqUofexppp}
 \U_g(t) \geq  \left( \frac{b}{f(\lambda t)}\right)^{\frac{1}{\alpha}}  =  \left( \frac{b}{\log(t) + o(\log(t))}\right)^{\frac{1}{\alpha}}.
\end{equation}
 \end{example}

\bigskip

It will also be useful to consider the inverse of the function $\U_g(t)$.  For $g(r) = ce^{br^{\alpha}}$ we claim that 
for any $0 < \alpha \leq 1$
\begin{equation}\label{InvUofexpr}
 u = \U_g(t)  \quad \Rightarrow \quad   t  \geq    \frac{c  u^{\alpha}}{\alpha  b}\cdot \exp\left({\frac{b}{u^{\alpha}}-1}\right),
\end{equation}
with  equality if $\alpha = 1$.   Indeed, we have 
$$
  u  =  \U_g(t) \geq   \frac{b}{f(\lambda t))}  \quad  \Leftrightarrow \quad 
 f(\lambda t) \geq \frac{b}{u^{\alpha}}.
$$
Using the relation $x ={\exp({f(x)})}/{f(x)}$, one gets 
$$
 \lambda t = \frac{\exp\left(f(\lambda t)\right)}{f(\lambda t)}  \geq  \frac{u^{\alpha}}{b} \exp \left( \frac{b}{u^{\alpha}}\right) , 
$$
which is equivalent to  \eqref{InvUofexpr} since  $\lambda = \frac{\alpha}{c}  \ee^{1/\alpha}$. 

\bigskip

\subsection{Asymptotic Behavior}
The $\U$-transform of a general function is generally not computable in a closed form,   but \eqref{solveU} gives us some useful estimates.
The first result in this direction is the following 
\begin{lemma}\label{lemmaineqU}
If $g(r)$  is  monotone increasing, then \  $\U_g(t) \leq \dfrac{1}{g^{-1}(t)}$. 
\end{lemma}

 \medskip 

\textbf{Proof.} 
Since $\tau(r)  \leq g(r)$ for any $r>0$,  we have 
$
  t = \tau(\rho(t)) \leq g(\rho(t))
$
for any $t$. Because $g$ is monotone increasing, so is its inverse, therefore  $g^{-1}(t) \leq \rho(t)$ for any $t$ and we conclude that 
$$
  \U_g(t) =  \frac{1}{\rho(t)}\left(1-\frac{t}{g(\rho(t))} \right) \leq  \frac{1}{\rho(t)} \leq  \frac{1}{g^{-1}(t)}.
$$
\qed

\medskip

\begin{lemma}\label{lemmalimitrho}
If $g(r)$  is unbounded and  monotone increasing, then $\displaystyle \lim_{t\to \infty} \rho(t) = \infty$.
\end{lemma}

\medskip

\textbf{Proof.}
Use again that $\tau(r)  \leq g(r)$ for any $r$. This implies that $t = \tau (\rho(t)) \leq g(\rho(t))$
and thus $\rho(t) \geq g^{-1}(t)$ and therefore $\displaystyle \lim_{t\to \infty} \rho(t) =  \lim_{t\to \infty} g^{-1}(t) = \infty$.

\qed

\bigskip 
 
The next result says that the $\U$-transform  of any function  $g: \R_+ \to \R_+$ satisfying \eqref{TSPG} (see Definition \ref{def.TSPG})
is asymptotically given by  ${1}/{g^{-1}(t)}$.
\begin{proposition} \label{prop.subexpg}
Let  $g: \R_+ \to \R_+$ be a strictly increasing differentiable function  satisfying \eqref{TSPG}, then  for  any $\varepsilon>0$, the following inequalities hold for $t$ large enough:
\begin{equation}\label{Uassy}
  \frac{1-\varepsilon}{g^{-1}(t)} \leq \U_g(t) \leq \frac{1}{g^{-1}(t)}.
\end{equation}
\end{proposition}

\bigskip

\textbf{Proof.}
The inequality $\U_g(t) \leq 1/g^{-1}(t)$ has already been proved in the previous Lemma. 
To prove the other inequality, we first  observe that \eqref{TSPG} has the following consequences: 
\begin{enumerate}[\quad (i)]
\item$\displaystyle \lim_{r\to \infty}   \frac{rg'(r)}{g(r)}= \infty$.
\item  $\displaystyle \lim_{r\to \infty} \frac{g(\lambda r)}{g(r)} = 0$ \  for any $\lambda \in [0,1)$.
\item  $\displaystyle \lim_{r\to \infty} \frac{g(\lambda r)}{\tau(r)} = 0$ \ for any $\lambda \in [0,1)$.
\end{enumerate}
Using (i) and the definition of $\tau(r)$,  we  obtain 
$$
 \lim_{t\to \infty}  \frac{t}{g(\rho(t))}   = \lim_{t\to \infty} \frac{\tau(\rho(t))}{g(\rho(t))} = \lim_{r\to \infty} \frac{\tau(r)}{g(r)} = 
\lim_{r\to \infty}  \frac{1}{1+ r\frac{g'(r)}{g(r)}}   =  0
$$ 
(we use Lemma \ref{lemmalimitrho} in the second equality). 
It  follows that for any ${\varepsilon}> 0$, one can find $t({\varepsilon})$ large enough so that 
$$
  \left(1 - \frac{t}{g(\rho(t))} \right)  \geq \left(1- \frac{\varepsilon}{2}\right)
$$ 
for any $t \geq t({\varepsilon})$. 

\smallskip 

On the other hand,  (iii) implies that for any $0 <\varepsilon < 1$,  there exists $r({\varepsilon})$ such that  \
$
 g((1-\varepsilon/2) r) \leq \tau(r) 
$ \  for any $r > r({\varepsilon})$.
Together with  the monotonicity of $g^{-1}$, this implies 
$$
   \left(1- \frac{\varepsilon}{2}\right) r \leq g^{-1} (\tau(r))
$$
 for any $r > r({\varepsilon})$. 
Enlarging $t({\varepsilon})$  if necessary, we may also assume that $\rho(t) \geq  r({\varepsilon})$,  for any $t \geq t({\varepsilon})$  (we use again Lemma \ref{lemmalimitrho} here), 
which then implies 
$$
    \left(1- \frac{\varepsilon}{2}\right) \rho(t) \leq  g^{-1}\left( \tau(\rho(t))\right) = g^{-1}(t),
$$
or equivalently
$$
  \left(1- \frac{\varepsilon}{2}\right) \frac{1}{g^{-1}(t)} \leq \frac{1}{\rho(t)}.
$$
We conclude that for any  $t \geq t({\varepsilon})$:
$$
  \U_g(t) =  \frac{1}{\rho(t)}\left(1-\frac{t}{g(\rho(t))} \right) \geq  \frac{1}{g^{-1}(t)}\left(1- \frac{\varepsilon}{2}\right)^2
\geq \frac{1}{g^{-1}(t)}\left(1-{\varepsilon}\right).
$$
\qed

\bigskip

\textit{Remark.}  Some of the above proofs are slightly delicate because  \textit{a priori} $\tau$ is not  invertible. The next lemma gives us a sufficient condition for this to be the case:

\begin{lemma}
 Let $g : [0, \infty) \to \R_+$ be a monotone function of class $C^2$ such $h(r) = \log (g(r))$ is concave, then the function 
$\tau$ is strictly increasing. 
\end{lemma}

 \textbf{Proof.} By concavity of  $h(r) = \log(g(r))$ we have 
$$
  \tau'(r) = \tau(r) \cdot  \frac{r (h'(r)^2-h''(r))}{1+rh'(r)} \geq  \tau(r) \cdot  \frac{r (h'(r)^2)}{1+rh'(r)} > 0.
$$
\qed

\subsection{Comparison with the Legendre Transform} \label{sec.legendre}

The $\U$-transform defined in   \eqref{defUT} is  equivalent to the Legendre transform, but is in some way more adapted to studying the asymptotic behavior of a function at infinity. Recall that the \textit{Legendre transform}  \index{Legendre transform}   (also called the \textit{Legendre--Fenchel transform}) of an arbitrary function $f : \R_+ \to \R$ is the function $ \LL_{f} : \R_+ \to \R_+ \cup \{+\infty\}$ defined as
$$
   \LL_{f}(y)  = \sup \left\{ yx - f(x) \mid  x\in \R_+ \right\}.
$$
The relationship between the $\U$-transform and the Legendre transform is the following  
\begin{lemma}
If  $f,g: \R_+ \to \R$ are two arbitrary  functions such that  $f(x)g \left(\frac{1}{x} \right) \equiv x$, then 
\begin{equation}\label{LUrelation}
  \U_{g}(t) = t \; \LL_{f}\left(\tfrac{1}{t}\right).
\end{equation}
\end{lemma}

The proof is elementary:  we have with $x= 1/r$ and $y=1/t$
$$
 \frac{1}{r}\left(1- \frac{t}{g(r)} \right) = t  \left(\frac{1}{tr} - \frac{1}{rg(r)}\right) 
= t \left(yx - f(x) \right).
$$
Taking the supremum in the above identity with respect to $r$ on the left hand side and $x$ on the  right hand side yields \eqref{LUrelation}.

\medskip

\paragraph{Acknowledgement.}  The authors are thankful to L. Bartholdi, J. Brieussel, N. Monod,  A.  Papadopoulos and C. Pittet for their comments and suggestions.

\end{document}